# DIOPHANTINE AND NON-DIOPHANTINE ARITHMETICS:

## Operations with Numbers in Science and Everyday Life

Mark Burgin

Department of Mathematics
University of California, Los Angeles
405 Hilgard Ave.
Los Angeles, CA 90095

> *[Arithmetic] is one of the oldest branches, perhaps the very oldest branch of human knowledge; and yet some of its most abstruse secrets lie close to its tritest truths.*
> H. J. S. Smith

**Abstract:** Science and mathematics help people better to understand world, eliminating different fallacies and misconceptions. One of such misconception is related to arithmetic, which is so important both for science and everyday life. People think that their counting is governed by the rules of the conventional arithmetic and that other kinds of arithmetic do not exist and cannot exist. It is demonstrated in this paper that this popular image of the situation with integer numbers is incorrect. In many situations, we have to utilize different rules of counting and operating. This is a consequence of the existing diversity in nature and society and to represent correctly this diversity people have to utilize different arithmetics. To distinct them, we call the conventional arithmetic Diophantine, while other arithmetics are called non-Diophantine. Theory of non-Diophantine arithmetics is developed in the book of the author "**Non-Diophantine arithmetics or *is it possible that* 2 + 2 *is not equal to* 4**." In this work, some properties of non-Diophantine arithmetics are considered, as well as their connections to numerical computations and contemporary physics are explained.

# 1. Introduction.

Sometimes theory in its development goes ahead of the boldest imagination. In this work, we are going to consider one of such situations when the entities that has been discovered and the theory that has been elaborated contrast basic dispositions of the common sense and contradict the acknowledged truth. However, observations and mathematical reasoning support these discoveries and theories, demonstrating their validity both in theoretical, practical and experimental aspects.

We begin with some historical examples, which demonstrate that bias and stereotypes are very powerful in society. History can help people to understand much better current situation and even to look into the future. Nevertheless, people very rarely take lessons from history. Hart (1971) even wrote a book "Why people don't learn from history."

So, here is one story from the history of cognition.

For thousands of years, people believed that the Earth is flat and suggested many hypotheses how it is situated. For example, one of them told that the Earth is situated on three elephants, three elephants are standing on the gigantic turtle, and this turtle is swimming in an infinite ocean. This is very beautiful image, but wrong in spite of all its beauty. At first philosophers suggested and then scientists proved that the Earth is round like a ball and rapidly moves in cosmic space.

Scientific misconceptions sometimes are even more powerful. For example, the Euclidean geometry was believed for 2200 years to be unique (both as an absolute truth and a necessary mode of human perception). People were not even able to imagine something different. In his "Metaphisical Foundations of Natural Science" (1786), the famous German philosopher Emmanuil Kant (1724-1804) argued that the world of science is a world of sense impressions arranged and controlled by the mind in accordance with innate categories such as space and time. Consequently, argued Kant, there never would be a way other than Euclidean geometry and Newtonian mechanics.

However, science and mathematics are functioning in a sequence of transformations. Sometimes these transformations change the whole foundation of human knowledge. In spite of this, almost unexpectedly some people began to understand that geometry is not unique.

Trying to improve the axiomatic system suggested for geometry by Euclid, three great mathematicians of the 19$^{th}$ century (Carl Friedrich Gauss (1777-1855), Nikolay Ivanovich Lobachewsky (1792-1856), and Janoš Bolyai (1802-1860)) discovered a lot of other geometries. At first, even the best mathematicians opposed this discovery and severely attacked Lobachewsky and Bolyai who published their results. Forecasting such antagonistic attitude, the first mathematician of his times Gauss was afraid to publish this. Nevertheless, progress of mathematics brought understanding and then recognition. This discovery is now considered as one of the highest achievements of the human genius. It changed to a great extent understanding of mathematics and improved comprehension of the whole world.

## 2. The Diophantine Arithmetic

In the 20$^{th}$ century, a similar situation existed in arithmetic of natural numbers, i.e. the numbers 1, 2, 3, ... . From all things that people are doing, counting is one of the most important. Without counting we cannot do a lot: we cannot develop science and technology, we cannot organize mass production, we cannot buy and sell, and so on and so forth. Every woman and every man, every boy and every girl perform counting many times a day. Calculators and computers were invented to help people to count. In the old days we used to say that computers do not do mathematics; they only do arithmetic. Later computer began to fulfill much more sophisticated tasks and to perform much more complex operations than we have in arithmetic. Some of its abilities look miraculous. However, counting lies at the bottom of all computer operations. Besides, for the purposes of some fields such as most molecular dynamics programs, computers only do arithmetic.

As writes Philip J. Davis (2000), "political, social, and economic applications of mathematics are now multiplying like Fibonacci's rabits." The majority of these applications involve numbers and arithmetic operations. To mention but one example, we can take the article "Water Arithmetic 'doesn't Adds" (Kirby, 2000), which deal with environmental problems. By standards of the research community of mathematicians, mathematics that people use in everyday life is trivial, "but all those everyday trivialities … made, consciously or unconsciously, philo- or miso-maths of most people."

People started using numbers and counting in prehistoric times. For thousands of years mathematicians studied numbers and counting learning a lot in this area. Now we know two mathematical disciplines studying numbers, relations between them and operations with them: arithmetic and number theory. Arithmetic takes the name from the Greek word for number (*arithmos*). However, the meaning of the word "arithmetic" changed through the ages. According to Smith (1963), some of the ancient philosophers and mathematicians used the term 'arithmetic' to denote the theory of numbers, which comprised such areas as the study of primes and properties of square numbers. To the practical use of numbers, including the methods of writing them, they gave the name logistic. This was looked upon by ancient scholars as somewhat a plebeian art, while arithmetic was considered an intellectual pursuit worthy of the attention of a philosopher. Other considered arithmetic as a whole, only having different aspects. For example, Plato in his Philebus writes: "Arithmetic is of two kinds, one is popular, and the other philosophical."

This situation is completely mirrored in the contemporary mathematics. Mathematicians assume that arithmetic, which comprises mostly representation and operation with numbers, is something completed long ago, trivial and aimed at the beginners, while number theory, which comprises investigation of properties of natural numbers, is a respected field of theoretical mathematics with its highly complicated,"deep" and abstract problems. However, knowing only one arithmetic of natural numbers, mathematicians forget that arithmetic (in the modern sense) is the base for theory of numbers. You change arithmetic and you need new theory of numbers as numbers change their properties. For example, the definition of a prime number is completely based on the definition of the operation of multiplication.

People's experience with numbers and, especially, with natural numbers is profound. All this time, much longer than that for the Euclidean geometry, people have believed and continue to believe that only one arithmetic of natural numbers exists.

However, do you think that we know everything about numbers and counting?

## 3. History of the Problem

Mathematical establishment treated the arithmetic of natural numbers as a primordial entity. For example, such prominent mathematician as Leopold Kronecker (1825-1891) wrote: "*God made the integers, all the rest is the work of man*".

Laymen have been even more persistent on this point of view. Almost all people had and have no doubts that **2 + 2 = 4** is the most evident truth in always and everywhere. Mathematicians supported this attitude. As it is written in The American Mathematical Monthly (April, 1999, p.375), "Although other sciences and philosophical theories change their 'facts' frequently, **2 + 2** remains **4**."

However, in our days, some people begin to doubt the absolute character of the ordinary arithmetic, where **2 + 2 = 4** and two times two is equal to four. As a result, scientists and mathematicians draw attention of the scientific and even general community to the foundational problems of natural numbers and the conventional arithmetic.

However, the roots of this fundamental problem of relevance of the conventional arithmetic may be found in ancient Greece. There was a group of philosophers, who were called Sophists and lived from the second half of the fifth century B.C.E. to the first half of the fourth century B.C.E.. Sophists asserted relativity of human knowledge and suggested many paradoxes, explicating complexity and diversity of real world. The famous philosopher Zeno of Elea ( 490-430 B.C.E.), who was said to be a self-taught country boy, invented very impressive paradoxes, in which he questioned the popular knowledge and intuition related to such fundamental essences as time, space, and number.

However, Greek sages posed questions, but in many cases, including arithmetic, suggested no answers. As a result, for more than two thousand years these problems were forgotten and everybody was satisfied with the conventional arithmetic. In spite of all problems and paradoxes, this arithmetic has remained very and very useful.

May be the first who questioned absolute validity of the conventional arithmetic in modern times was the famous German scientist Herman Ludwig Ferdinand von Helmholtz (1821-1894). In his "Counting and Measuring" (1887), Helmholtz considered an important problem of applicability of arithmetic, at that time it was only one arithmetic, to physical phenomena. This was a natural approach of a scientist, who even mathematics judged by the main criterion of science – observation and experiment.

His first observation was that as the concept of number is derived from some practice, usual arithmetic is applicable to these experiences. However, it is easy to find many situations when this is not true. To mention but a few described by Helmholtz, one raindrop added to another raindrop does not make two raindrops. In a similar way, when one mixes two equal volumes of water, one at 40° Fahrenheit and the other at 50°, one does not get two volumes at 90°. Alike the conventional arithmetic fails to describe correctly the result of combining gases or liquids by volumes. For example (Kline, 1980), one quart of alcohol and one quart of water yield about 1.8 quarts of vodka.;

Later the famous French mathematician Henri Lebesgue facetiously pointed out (cf. Kline, 1980) that if one puts a lion and a rabbit in a cage, one will not find two animals in the cage later on.

However, very few (if any) paid attention to the work of Helmholtz on arithmetic, and as still no alternative to the conventional arithmetic has been suggested, these problems were once more forgotten.

It took almost a hundred of years to revive them in our times. The most extreme view is that there is only a finite quantity of natural numbers. It is one of the central postulates of ultraintuitionism (Yesenin-Volpin, 1960). Van Danzig in his article (1956) explained why only some of natural numbers may be considered finite. Consequently, all other mathematical entities that are called traditionally natural numbers are only some expressions but not numbers. These arguments are supported and extended in (Blehman *et al.*, 1983). In addition this, we know that any computer arithmetic contains only a finite set on numbers.

Other authors are more moderate in their criticism of the conventional arithmetic. They write that not all natural numbers are similar in contrast to the presupposition of the conventional arithmetic that the set of natural numbers is uniform (Kolmogorov, 1961; Littlewood, 1953; Birkhoff and Barti, 1970; Rashevsky, 1973; Dummett 1975; Knuth, 1976). Different types of natural numbers have been introduced, but without changing the conventional arithmetic. For example, one of the greatest mathematicians of the 20$^{th}$ century

Andrei Nikolayevich Kolmogorov (1961) suggested that in solving practical problems it is worth to separate *small, medium, large,* and *super-large* numbers.

A number *A* is called *small* if it is possible in practice to go through and work with all combinations and systems that are built from *A* elements each of which has two inlets and two outlets.

A number *B* is called *medium* if it is possible to go through down and work with this number. However, it is impossible to go through and work with all combinations and systems that are built from *B* elements each of which has two or more inlets and two or more outlets.

A number *C* is called *large* if it is impossible to go through a set with this number of elements. However, it is possible to elaborate a system of denotations for these elements.

If even this is impossible, then a number is called *super-large.*

According to this classification, 3, 4, and 5 are small numbers, 100, 120, and 200 are medium numbers, while an example of a large number is given by the quantity of all visible stars. Really, is we invite 4 people, we can consider all their possible positions at a dinner table. If you come to some place where there are 100 people, you can shake hands with everybody. Although, it might take too much time. What concerns the visible stars, you cannot count them, although, a catalog of such stars exists. Using this catalog, it is possible to find information about any of these stars.

This classification of numbers is based on our counting abilities. Consequently, borders between classes are vague and unstable. Higher counting abilities make borders between classes higher. For example, 10 is a medium number for an ordinary individual, but a small number for a computer. However, some numbers belong to a definite class of this typology in all known situations. For example, 300 is a medium number both for people and computers.

In a similar way to what has been done by Kolmogorov and on akin grounds, the outstanding British mathematician John Edensor Littlewood (1953) separated all natural numbers into an infinite hierarchy of classes.

All mathematicians who were wise enough to distrust the complete adequacy of the conventional arithmetic may be divided into two groups. Some write (like Helmholtz, Kolmogorov or Littlewood) that on practice natural numbers and operations with them are different from those that we know from mathematics. Others admit (like Kline, 1967) that different arithmetics exist but we do not know what they are.

The well-known Russian mathematician Pyotr Konstatinovich Rashevsky (1973) was one of the firsts to formulate explicitly the construction problem for a theory of an arithmetic of natural numbers that is different from the conventional one. He even made some predictions about such arithmetics:

*"It must not be expected that* [such] *hypothetical theory, if it would be ever destined to se the light of day, will be unique; on the contrary, it will have to depend on certain "parameters" (with a role distantly reminiscent of the radius of Lobachevsky space when we repudiate Euclidean geometry in favor of non-Euclidean). It may be expected that in the limiting case the hypothetical theory should coincide with the existing one."*

More recently, the American mathematician Brian Rotman (1997) also directly formulated the problem of elaboration of arithmetics that were essentially different from the conventional one. He based his suggestions on a series of examples demonstrating that many laws of the conventional arithmetic are not true in different situations. Rotman calls new structures non-Euclidean arithmetics, although he does not describe them. However, it is more natural to call the conventional arithmetic by the name *'Diophantine arithmetic'* than *'Euclidean arithmetic'* because Diophantus contributed much more to the development of arithmetic than Euclid. Consequently, new structures acquired the name "*non-Diophantine arithmetics*'.

The main arguments of Morris Kline concerning arithmetic are similar to those from (Rashevsky, 1973; Burgin, 1977). Although, the examples supporting these arguments are different. Let us consider one of these examples, which is taken from the book of Morris Kline (1967).

If a farmer has two herds consisting of 10 and 25 heads of cows, respectively, he knows by adding 10 and 25 that the total number of cows is 35. That is, he need not count his cows. Suppose, however, he brings the two herds of cows to market where they are selling for $100 apiece. Will a herd of 10 cows which might bring $1000 and a herd of 25 cows which might bring $2500 together bring in $3500? Every businessman knows that when supply exceeds demand, the price may drop, and hence 35 cows may bring in only $3000. In some idealized world the value of the cows may continue to be $3500, but in actual situations this need not be true.

Consequently, continues Kline (1967), mathematicians are, of course, free to introduce the symbols 1, 2, 3, … , where 2 means 1+1, 3 means 2+1, and so on. We can even deduce from this that 2 + 2 = 4. But the question is not whether the mathematician can set up definitions and axioms and deduce conclusions. It is necessary to know whether this system necessarily expresses truths about the physical world.

According to Kline (1967), discovery of non-Euclidean geometries had taught mathematicians that geometry does not offer ultimate truths. That was the reason why many turned to the ordinary number system and the developments built upon it and maintained that this part of mathematics still offers unquestionable truths. The same thought is often expressed today by people who, wishing to give an example of an absolute truth, quote 2 + 2 = 4. However, examination of the relationship between our ordinary number system and the physical situations to which it is applied vividly demonstrates that it does not offer truths.

**4. Discovery of Non-Diophantine Arithmetics**

The conventional arithmetic may be called Diophantine because the ancient Greek mathematician Diophantus was the first who made an essential contribution to arithmetic. Diophantus, often called the 'father of algebra', is best known for his *Arithmetica*, a treatise on the theory of numbers. Diophantus did his work in the great city of Alexandria. At this time Alexandria was the center of mathematical learning. However, very little is known of Diophantus's life and there has been much debate regarding the date at which he lived. According to some sources, he lived 84 years and had one son.

There are a few limits that can be put on the dates of Diophantus's life by analyzing his references. On the one hand Diophantus quotes the definition of a polygonal number from the work of Hypsicles so he must have written this later than 150 B.C.E. On the other hand Theon of Alexandria, the father of Hypatia, quotes one of Diophantus's definitions so this means that Diophantus wrote no later than 350 C.E. However this leaves a span of 500 years, so we have not narrowed down Diophantus's dates a great deal by these pieces of information. The reason is that while many references to the work of Diophantus have been made, Diophantus himself, made very few references towards other mathematicians' work, thus making the process of pinpointing his dates harder.

In the same way as it was with the Euclidean geometry, the Diophantine arithmetic has been unique and nonchallengable for a very long time – people have not known other arithmetics. Its position in human society has been and is now even more stable and firm than the position of the Euclidean geometry before the discovery of the non-Euclidean geometries. Really, all people use the Diophantine arithmetic for counting. Just arithmetical operations do out of all people some kind of consumers of mathematics. At the same time, Euclidean geometry is only studied at school and in real life rather few specialists use it. It is arithmetic, and not geometry, which is considered as a base for the whole mathematics in the intuitionistic approach.

By the well-known English mathematician Henry John Stanley Smith (1826-1883), arithmetics (and namely, the Diophantine arithmetic) is one of the oldest branches, perhaps the very oldest branch, of human knowledge. His older contemporary, very talented German mathematician Carl Gustav Jacob Jacobi (1805-1851) said: "*God ever arithmetizes*".

Nevertheless, in spite of such a high estimation of the Diophantine arithmetic, its uniqueness and indisputable authority has been recently challenged. The first family of non-Diophantine arithmetics was discovered in 1975 (Burgin, 1977). Like geometries of Lobachewsky, these systems, which are called *projective arithmetics*, depend on a special parameter. Although, this parameter is not a number as in the case of Lobachewsky geometries. Arithmetics have a functional parameter, that is, any arithmetic in this family, properties and laws of its operations depend on a definite function *f(x)*. The Diophantine arithmetic is a member of this family: its parameter is equal to the identity function $f(x) = x$.

Later another family of non-Diophantine arithmetics, which are called *dual arithmetics*, was introduced (Burgin, 1980). This family also has a functional parameter, and the Diophantine arithmetic is also a member of this family: its parameter is equal to the identity function $f(x) = x$. However, many properties of the arithmetics from the second family are essentially different in comparison with arithmetics from the first family. For example, in projective arithmetics, we can have an identity $n + 1 = n$. In dual arithmetics, $n + 1$ is always larger than $n$. The book "**Non-Diophantine arithmetics or *is it possible that* 2 + 2 *is not equal to* 4**" (Burgin, 1997) contains a detailed study of non-Diophantine arithmetics from both families.

### 5. Properties of Non-Diophantine Arithmetics

In some non-Diophantine arithmetics, even the most evident truth (such as **2 + 2 = 4** or two times two is equal to four) may be discarded. Some of them (projective arithmetics) possess similar properties to those of transfinite numbers arithmetics built by the great German mathematician Georg Cantor (1845-1918). For example, a non-Diophantine arithmetic may have a sequence of numbers $a_1, a_2, ..., a_n, ...$ such that for any number $b$ that is less than some $a_n$ the equality $a_n + b = a_n$ is valid. This is an important property of infinity, which is formalized by transfinite (cardinal and ordinal numbers). The equality $a^2 = a$ is another interesting property of some transfinite numbers. Such equality may be also true in some non-Diophantine arithmetics. Thus, non-Diophantine arithmetics provide mathematical models in which finite objects – natural numbers – acquire features of infinite objects – transfinite numbers. In such a way, it is possible to model and to describe behavior of infinite entities in finite domains.

In non-Diophantine arithmetics, it is also possible that two times two is not equal to four. We can add different numbers to the same number and get the same result. There are such non-Diophantine arithmetics that have the largest number. All this contrasts what is possible in the Diophantine arithmetic.

Let us consider some other features of non-Diophantine arithmetics and compare them with the properties of the Diophantine arithmetic.

We know from school that the main laws of the Diophantine arithmetic are:

1. Commutativity of addition: $a + b = b + a$;

2. Associativity of addition: $(a + b) + c = a + (b + c)$;
3. Commutativity of multiplication: $a \cdot b = b \cdot a$;
4. Associativity of multiplication: $(a \cdot b) \cdot c = a \cdot (b \cdot c)$;
5. Distributivity of multiplication with respect to addition: $a \cdot (b + c) = a \cdot b + a \cdot c$;
6. Zero is a neutral element with respect to addition: $a + 0 = 0 + a = a$;
7. One is a neutral element with respect to multiplication: $a \cdot 1 = 1 \cdot a = a$.

Now we may ask whether these laws are valid for non-Diophantine arithmetics. First, addition and multiplication are always commutative, zero is a neutral element with respect to addition, and one is a neutral element with respect to multiplication in all non-Diophantine arithmetics. At the same time, the laws of associativity and distributivity fail in the majority of non-Diophantine arithmetics. Only special conditions on the functional parameter of the non-Diophantine arithmetic in question provide validity of these laws. These conditions are obtained in (Burgin, 1997).

Second, the Diophantine arithmetic possesses the, so-called, Archimedean property, which is important for proofs of many results in arithmetic and number theory. It states that if we take any two natural numbers *m* and *n*, in spite that *n* may be enormously larger than *m*, it is always possible to add *m* enough times to itself, i.e., to take a sum *m* + *m* + … + *m*, so that the result will be larger than *n*. This property is also invalid in the majority of non-Diophantine arithmetics. As we have seen above (cf. the example "*A Woman and a Pay Phone*"), it is possible to add 1 to itself as many times as you can but never get 5. The Archimedean property is important for proving that the set of all natural numbers is infinite as well as the set of all prime numbers. Thus, having in general no Archimedean property in non-Diophantine arithmetics, we encounter such arithmetics that have only finite number of elements or such infinite arithmetics that have only finite set of prime numbers (Burgin, 1997).

Another example of a non- Archimedean arithmetic in the nonstandard arithmetic (Robinson, 1966). However, this arithmetic includes infinitely big and infinitely small numbers, going thus, beyond natural numbers.

Third, the Diophantine arithmetic always contains infinitely many numbers. Some of non-Diophantine arithmetics contain the largest number and thus, are finite. However, this is possible only for perspective arithmetics. All dual arithmetics are infinite. So, we can see that many properties of non-Diophantine arithmetics contradict our intuition and contrast to what we know from our experience.

**6. When and Where We Need Non-Diophantine Arithmetics**

However, the peculiarity of the situation is that, as we will see, in a lot of cases people unconsciously use non-Diophantine arithmetics. In spite of this, it is so hard for people to understand and what is even more difficult to accept non-Diophantine arithmetics. Power of people's stereotypes is vividly demonstrated by the book (Blehman *et al.*, 1983). At first (section 1.2.4), the authors of the book explain with many examples and references that our intuition of natural numbers and arithmetic are very misleading. After this (section 1.2.5), they announce that it is completely impossible that two times two is not equal to four. The authors are even trying to prove this utilizing the following probabilistic reasoning. Here are their arguments (p.50).

*Really, the statement that two times two is equal to four may be taken as an example of the most evident truth. Although, nobody doubts that this is a true equality, it is possible to evaluate formally probability that in reality two times two is equal to five, while the standard statement that two times two is equal to four is a result of a constantly repeated arithmetical mistake. Let us suppose that any individual performing multiplication with numbers that are less than ten can decrease the result by one with the probability $10^{-6}$. This corresponds to several such mistakes during his or her life. If we assume that that through the whole history of mankind, $10^{10}$ people performed the multiplication "two times two" $10^{6}$ times during the life of each of them, then the probability that they repeated this mistake of decreasing the result is less than $10^{-10^{17}}$. Thus,* the authors conclude, *the probability is so small that the event is absolutely impossible and we see that two times two is equal to four.*

However, there are such non-Diophantine arithmetics in which two times two is not equal to four. Besides, the authors are absolutely sure that nobody ever did the mistake. Both these facts show to what extent probabilistic proofs may be misleading. This reminds us how some people try to prove that it is impossible that life emerged by itself. Actually, they only calculate some small probability of this event and find that this probability is very small. Then, in spite that some of their assumptions are rather dubious, they claim that all living organisms were created.

When we consider non-Diophantine arithmetics, it is possible to think that they are absolutely formal constructions like many other mathematical objects, which are very far from the real world. But let us recollect that a similar skepticism and mistrust met the discovery of non-Euclidean geometries. Even Gauss (in spite of being acknowledged as the greatest mathematician of his time) did not dare to publish his results concerning these geometries as he was not able to find anything that is similar to them in nature. Lobachewsky called his geometry an imaginable one. But afterwards it was discovered that the real physical space fits non-Euclidean geometries, and that the Euclidean geometries do not have such essential applications as the non-Euclidean ones. In this aspect, the situation with non-Diophantine arithmetics is different. In spite of the short time, which has passed after their discovery, it has been demonstrated that many real phenomena and processes exist that match the non-Diophantine arithmetics.

As a matter of fact, much earlier than non-Diophantine arithmetics appeared, Littlewood (1953) considered an example demonstrating how the rules of non-Diophantine arithmetics (in spite of that they were unknown at that time) can be imposed upon the real world. Several similar and even more lucid examples are exposed in (Davis and Hersh, 1986) and in (Kline, 1967). Let us consider some of such situations.

1. A market sells a can of tuna fish for $1.05 and two cans for $2.00. So, we have $a + a \neq 2a$.

2. In a similar way, coming to a supermarket, you may bye one gallon of milk for $2.90 while two gallons of the same milk will cost you only $4.40. Once more, we have $a + a \neq 2a$.

3. Even more, coming to a supermarket, you can see an advertisement "Buy one, get one free." It actually means that you can buy two items for the price of one. Such advertisement may refer almost to any product: bread, milk, juice etc. For example, if one gallon of orange juice costs $5, then we to the equality 5 + 5 = 5. It is impossible in the conventional arithmetic but it true for some non-Diophantine arithmetics.

4. To make the situation, when ordinary addition is inappropriate, more explicit, an absurd but not unrelated question is formulated: If the Mona Lisa painting is valued at $10,000,000, what would be the value of two Mona Lisa paintings.

5. Another example: when a cup of milk is added to a cup of popcorn then only one cup of mixture will result because the cup of popcorn will very nearly absorb a whole cup of milk without spillage. So, in the last case we have **1 + 1 = 1**. It is impossible in the conventional arithmetic but it true for some non-Diophantine arithmetics.

Sometimes, we encounter a situation such that we cannot use the Diophantine arithmetic. As an example, we consider the situation, which we call "A Woman and a Pay Phone." A woman going to make a phone call comes to a pay phone. She has enough money – ten dollars, but all her money are in one cent coins. Coming to the phone closer, the woman reads that the phone accepts only nickels, dimes, and quarters and it is necessary to have 35 cents to make a local call. So, adding his cents, the woman is not able to make neither necessary 35 cents, nor a quarter, nor a dime, nor even a nickel. In other words, in the arithmetic of this bounded but real world (especially, if there are no other people close to that place), adding 1 to 1 to 1 and so on will never give 5. It is necessary to go outside this small world, to ask other people or to use an automaton in order to change separate cents for nickels, dimes, and quarters. But then it will be a different arithmetic.

In an appropriate non-Diophantine arithmetic, this means that 5 is much bigger than 1 (5 >> 1), that is, by adding 1 to 1 to 1 and so on, it is impossible to get 5. In the Diophantine arithmetic, any number is a sum of 1 plus 1 plus 1 and so on.

As a rule different automata change the Diophantine arithmetic. For example, a computer arithmetic is a non-Diophantine arithmetics. This is a result of round-off procedures and existence in such an arithmetic the largest number. Consequently, if we want to build better models for numerical computations than we have now, it is necessary to utilize in such models relevant non-Diophantine arithmetics.

These examples demonstrate that non-Diophantine arithmetics are important for business and economics. Some economical problems and inconsistencies caused by the conventional arithmetic are considered in (Tolpygo, 1997). As some studies of economy show, some finite quantities possess properties of infinite numbers with respect to people's practice (cf., also, Birkhoff and Barti, 1970). Consequently, when one applies mathematics to solve such problems, the results are often mathematically correct but practically misleading. Utilization of non-Diophantine arithmetics eliminates those problems and inconsistencies.

In addition to this, non-Diophantine arithmetics solve some problems that remained unsolved from the time of ancient Greece. One of such problems is discussed in (Rashevsky, 1973). This is the so-called, "paradox of a heap," which is attributed to Zeno.

Let us consider a heap of grains. If we add to this heap one grain, the heap is not changing. Consequently, if we take the number $K$ of the grains in this heap, then adding 1 to $K$ does not change $K$. This contradicts the main law of the Diophantine arithmetic stating that for an arbitrary number $k$, $k + 1$ is not equal to $k$, and gives birth to a paradox if we have only one arithmetic. Non-Diophantine arithmetics solve this paradox.

A reader may ask what for we are trying to do something with puzzles that were suggested thousands years ago and look artificial to a modern reader. However, the paradox of a heap has a direct analogy in our times both in science and everyday life. For example, you are buying a car for $30,000. Then suddenly, when you have to pay, the price is changed and becomes one cent greater. Do you think that the new price is different from the initial one or you consider it practically one and the same price? It is natural to suppose that any sound person has the second opinion. Consequently, we come to the same paradox: if $k$ is the price of the car in cents, then in the Diophantine arithmetic $k + 1$ is not equal to $k$, while in reality they are the same. Moreover, imagine that you are going to receive your salary in cents, the sum of which will be equal to the amount that you receive now but which will be given to you once a year. Do you think it will be the same salary or not? For many people, there is a

difference how often they receive their salary, although the sum remains the same. As it is possible to show, non-Diophantine arithmetics solve these paradoxes.

There are situations when very old problems and paradoxes appear on a new level of cognition. This happened with many profound ideas of ancient philosophers. All educated people know that material things are built of atoms. However, idea of atoms was introduced much earlier than atoms were really discovered in nature. Outstanding philosophers Democritus and Levkipus from ancient Greece elaborated idea of atoms as the least particles of matter. For a long time this idea was considered false due to the fact that scientists were not able to go sufficiently deep into the matter. Nevertheless, the development of scientific instruments and experimental methods made possible to discover such micro-particles that were and are called atoms although they possess very few of those properties that were ascribed to them by ancient philosophers.

Another great idea of ancient Greece was the world of ideas existence of which was postulated by Plato. In spite of the attractive character of this idea, the majority of scientists and philosophers believe that the world of ideas does not exist, because nobody had any positive evidence in support of it. The crucial argument of physicists is that the main methods of verification in modern science are observations and experiments, and nobody has been able to find this world by means of observations and experiments. Nevertheless, some modern thinkers, including such outstanding sages as philosopher Karl Popper, mathematician Kurt Gödel, and physicist Roger Penrose, continued to believe in the world of ideas giving different interpretations of this world but suggesting no ways for their experimental validation.

However, science is developing, and this development provided recently for the discovery of the world of structures (Burgin, 1996). On the level of ideas, this world may be associated with the Platonic world of ideas in the same way as atoms of the modern physics may be related to the atoms of Democritus. Existence of the world of structures is proved by means of observations and experiments (Burgin, 1997a). This world of structures constitutes the structural level of the world as whole. Each system, phenomenon, or process either in nature or in society has some structure. These structures exist like such things as tables, chairs, or buildings do and form the structural level of the world. When it is necessary to investigate or to create some system or process, it's possible to do this only by means of knowledge of the corresponding structure. Structures determine the essence of things in the same way as Aristotle ascribed to forms of things. Consequently, structures unite ideas of Plato with forms of Aristotle, eliminating contradictions that existed between these teachings.

In a similar way, paradoxes from ancient Greece are sometimes revived in modern physics. For example, theory of chaos encounters now many difficult problems caused by insufficiency of modern mathematics to represent correctly chaotic behavior. A chaotic behavior is defined as such behavior in which arbitrary small changes in external parameters or/and initial conditions cause essential changes in the behavior of a system. Such definition, as the expert in the theory of chaos Vladimir Gontar writes (Gontar, 1993; Gontar and Ilin, 1991), makes it problematic to use the apparatus of differential calculus to describe chaotic motion mathematically. The main obstacle is that solutions of differential equations must be continuous and robust. This condition implies that small changes in the initial conditions of any parameter may result only in small changes of trajectories. Utilization of numerical computations for solving problems of the system dynamics description cast doubt on the robustness of iterative methods and raised the question of the accuracy of numerical computations. This problem leads us, as emphasizes Vladimir Gontar (1993), to the paradoxes of ancient Greeks: is it mathematically and logically possible to formulate a contradiction-free description of the process of approaching an object when the distance to this object contains an infinite number of segments, involving a infinite number of steps necessary to reach this object?

Non-Diophantine arithmetics suggest a new understanding for this problem. It is possible (and Heisenberg's principle of uncertainty supports it) that we can come only to a definite distance to the object in question. It is possible to do this in a finite number of steps. All consequent steps cannot make the distance to the object smaller.

## 7. More Properties of Non-Diophantine Arithmetics

An interesting peculiarity of the new arithmetics is their ability to provide means for mathematical grounding of some intuitive constructions used by physicists. As writes Rashevsky (1973), "*since physicist uses only the apparatus proffered to him by the mathematician, the absolute power of natural numbers spreads to physics, and by means of the real axis predetermines to a considerable extent the possibilities of physical theories.*" As an example of new possibilities provided by non-Diophantine arithmetics, we can consider a mathematical model for such binary relations as "much larger" (denoted by >> ) and "much lesser" (denoted by << ). These relations are very frequently used in physics but there they have no exact meaning without non-Diophantine arithmetics. At the same time, they have a natural formalization in non-Diophantine arithmetics.

In a non-Diophantine arithmetic, we have $a << b$ if $b + a = b$. This relation is closely connected to the Archimedean property of the Diophantine arithmetic, which has been considered above. In particular, the following result is true.

**Theorem.** `A non-Diophantine arithmetic` **A** `has the Archimedean property if and only if the relation` $a << b$ `in` **A** `always implies` $a = 0$.

However, this relation may be true only in projective arithmetics. All dual arithmetics are Archimedean. In particular, we have $b + a > b$ for any two elements $a > 0$ and $b$ from an arbitrary dual arithmetic *A*. This implies that any dual arithmetic is infinite and we can obtain any number in such an arithmetic by adding 1 to 1 to 1 and so on. At the same time, dual arithmetics possess a binary relation ‹‹ , which is similar to the relation << of the projective arithmetics. Namely, we have $a ‹‹ b$ if $b + a = b + 1$. Looking very similar, these relations are, nevertheless, intrinsically different: while the relation $a << b$ implies that addition of *a* to *b* does not change *b* at all, the relation $a ‹‹ b$, means that addition of *a* to *b* changes *b* as little as possible.

Another interesting relation in projective arithmetics is $<<<$. Namely, we have $a <<< b$ if $b \circ a = b$ where $\circ$ denotes operation of multiplication in the corresponding arithmetic. As it is easy to see, $1 <<< b$ for any element $b$.

Other peculiarities of physical theories become clearer in the light of non-Diophantine arithmetics. For example, in physics there are some absolute values like the speed of light C or the absolute zero temperature. According to relativity theory, C is the largest speed attainable by material bodies. If you add some speed V to C, you have the same C, in other words, $C + V = C$. This equality is possible in some non-Diophantine arithmetics but impossible in the conventional, Diophantine arithmetic. In terms of non-Diophantine arithmetics, this means that C is much larger than any number, even than C itself. In this peculiarity of relativity theory, we can clearly see the old paradox of a heap, which is considered above.

In addition to this, finiteness of all material bodies in the universe suggests that temperature has not only the lower bound but also the upper bound, and we again come to non-Diophantine arithmetics.

That is why, some physicists (cf., for example, Zeldovich *at al*, 1990) emphasized that fundamental problems of modern physics are dependent on our ways of counting. This idea correlates with problems of modern physical theories in which physical systems are described by chaotic processes. Taking into account the fact that chaotic solutions are obtained by computations, physicists ask (Cartwrite and Piro, 1992; Gontar, 1997) whether chaotic solutions of the differential equations, which model different physical systems, reflect the dynamic laws of nature represented by these equations or whether they are solely the result of an extreme sensitivity of these solutions to numerical procedures and computational errors.

Other examples of situations for which non-Diophantine arithmetics give more adequate models than the Diophantine arithmetic are related to the concept of convergence for series and consequently, to computations and numerical analysis. As it was emphasized by the great French mathematician Henri Poincare (cf. Blehman *et al*., 1983), convergence for a physicist and a mathematician does not coincide. Really, let us consider a simple example of two series: $1000^n/n!$ and $n!/1000^n$. Mathematicians will call the first series convergent and the second series divergent because its member can grow without limits. At the same time, astronomers will say that the first series is divergent because its first 1000 member are increasing, while the second series is convergent because its first 1000 member are decreasing and rather fast.

Poincare concludes that both opinions are legal: the first one in theoretical studies and the second one in practical applications. It is important to separate these areas.

In other words, a series that converges in the sense of mathematical analysis may be divergent from the point of view of an astronomer and vice versa. For example, if some series converges but so slowly that to observe this convergence it is necessary to make more computational operations than the fastest computer can do in ten years. Then such series is practically divergent. At the same time, divergent series, which are called asymptotic series, are successfully used in physics (cf., for example, Collins, 1984).

However, if we have only one arithmetic, then the concept of convergence has to be unique by its definition. When, on the contrary, we have different arithmetics, then one and same definition of convergence may give unalike results. The reason for these differences is that series and their convergence are considered in different arithmetics. It is possible to show that arithmetics that are used for applications in many fields are, as a rule, non-Diophantine.

## 8. How We Accept New Ideas

However, in spite of all evident facts that prove existence and usage of non-Diophantine arithmetics, people's conservatism, inertia, and bias prevent them from empowering themselves with this tool of cognition and practical activity. An important regularity of society is that its life is based on a multiple stereotypes. This is also true for science (Burgin and Onoprienko, 1996). There are many historical examples of such situations.

One of the most notorious is from the history of the discovery of non-Euclidean geometries, which was mentioned before. Not only laymen did not want even to listen about this outstanding achievement of human mind, but even well-known and authoritative mathematicians opposed this new theory. For example, Ostrogradsky, who was a contemporary of Lobachevsky and was considered the best Russian mathematician at that time, wrote a withering review on the work of Lobachevsky on non-Euclidean geometries. In addition to this, a pamphlet on Lobachevsky was published in the most important Russian newspaper of that time. The main thesis of this pamphlet was that some people come from the province, call themselves mathematicians, and write very clumsy and awkward texts about such things that are impossible because they are never possible.

We know that Kant claimed that (Euclidean) geometry is given to people apriory, i.e., without special learning. However, one might say that Kant was a philosopher and not an expert in mathematics. Besides, when Kant wrote about geometry, non-Euclidean geometries were not yet built. Only some people suggested hypotheses that other geometries might exist. However, William Rowan Hamilton (1805-1865), certainly one of the outstanding mathematicians of the 19$^{th}$ century, expressed similar consideration in 1837 when the works of Lobachevsky (1829) and Boyai (1932) had been published. Hamilton said:

"No candid and intelligent person can doubt the truth of the chief properties of *Parallel Lines*, as set forth by Euclid in his *Elements*, two thousand years ago; though he may well desire to see them treated in a clearer and better method. The doctrine involves no obscurity nor confusion of thought, and leaves in the mind no reasonable ground for doubt, although ingenuity may usefully be exercised in improving the plan of the argument."

Even in 1883, another outstanding mathematician Arthur Cayley (1821-1895) in his presidential address to the British Association for the Advancement of Science affirmed:

"My own view is that Euclid's twelfth axiom [usually called the fifth or parallel axiom] in Pfayfair's form of it does not need demonstration, but is part of our notion of space, of the physical space of our experience…"

Moreover, it is known (Kline, 1980) that even in 1869 and 1870 the noted French mathematician Joseph Bertran (1822-1900) published in Comptes Rendus (Bulletin) of the Paris Academy of Sciences his works in which he disproved non-Euclidean geometries.

Outstanding achievements of the great French mathematician Évariste Galois (1811-1832) give another example of such situation. As writes E.T. Bell (1965), "in all the history of science there is no completer example of the triumph of crass stupidity over untamable genius than is afforded by the too brief life of Évariste Galois."

The first time Galois submitted his fundamental discoveries to the French Academy of Sciences, Cauchy, who was the leading French mathematician of that time, promised to present this, but he forgot and even lost the abstract of the works. The second time Galois submitted his epochal papers on the theory of algebraic equations. These ideas later gave birth to the modern algebra, in which structures were studied instead of equations. However, these works were also lost because nobody understand either the essence of Galois discoveries or even their importance. Only by a mere chance his works were published long after his death.

Still another example of the stereotypes that hinder acceptance of highly original and important theories and ideas is the attitude to the works of the outstanding German logician Gotlob Frege (1848-1925). Grattan-Guinness writes in his History of the Mathematical Sciences (1998):

"*Neglect* [of Frege] *during his lifetime is normally explained as due to his unusual notation ... but I learnt it in 20 minutes and I cannot see why his contemporaries would have needed longer.*"

The life and works of such great English physicist as Oliver Heaviside (1850-1925) give one more example and we can find a lot more.

What concerns non-Diophantine arithmetics, here are we describe only two situations. The first one is related to one of the well-known contemporary physicists. When somebody came to him with dubious results in theoretical physics, the physicist said that that results are not possible and continued:

"*You might as well claim that you can prove that* $2 + 2 = 3$. *If you think you can do that, don't waste my time, go to the market and buy two* $2 *items and try to pay for them with* $3."

Some time after this episode, this physicist was informed about non-Diophantine arithmetics, in which it is possible that $2 + 2 = 3$, and given an example when it was possible " to go to the market and buy two $2 items and try to pay for them with $3." This gives an experimental proof that it *is possible* that $2 + 2 = 3$.

The reaction of the physicist was very characteristic. Although he consented that it is possible to "*go to the market, buy two* $2 *items, and pay for them with* $3," he rejected the possibility of other arithmetics. Here is exactly what he responded:

"*I just checked on my fingers*: 2 + 2 *still* = 4. *Therefore, I have proven experimentally that the theory* [of non-Diophantine arithmetics] *is wrong.*"

Then he added, "*Mathematicians are fond of calling an apple a banana.*"

At first, let us look at his "as-if proof." In a same way, some people "prove" that there are neither atoms nor subatomic particles. They say, "I just looked around and did not see any atoms. Therefore, I have proven experimentally that the theory [of atoms] is wrong."

Then such a hypothetical person may add, "Physicists are fond of inventing things that do not exist."

A similar situation with "as-if proofs" is described in the following joke, which is aversion of the one from (Polya, 1956). The joke "How a mathematician, physicist, engineer, and computer scientist prove that all odd integers are prime."

The mathematician says, " Well, 1 is prime, 3 is prime, 5 is prime, and by induction, we have that all the odd integers are prime."

The physicist then says, "I'm not sure of the validity of your proof, but I think I'll try to prove it by experiment." He continues, "Well, 1 is prime, 3 is prime, 5 is prime, 7 is prime, 9 is ... uh, 9 is an experimental error, 11 is prime, 13 is prime... Well, it seems that you're right."

The engineer then says, "Well, actually, I'm not sure of your answer either. Let's see... 1 is prime, 3 is prime, 5 is prime, 7 is prime, 9 is ..., 9 is ..., well if you approximate, 9 is prime, 11 is prime, 13 is prime... Well, it does seem right."

Not to be outdone, the computer scientist comes along and says "Well, you two have got the right idea, but you have end up taking too long doing it. I've just whipped up a program to REALLY go and prove it..." He goes over to his terminal and runs his program. Reading the output on the screen he says, "3 is prime, 5 is prime, 7 is prime, 7 is prime, 7 is prime, 7 is prime, 7 is prime...."

In another version of the joke, the physicist then says, "I will prove that all odd integers are prime by experiment. Well, 1 is prime, 3 is prime, 5 is prime, 7 is prime, 9 is ... uh, 9 is equal 3 times 3 and as 3 is prime, we may, as well, call 9 prime, 11 is prime, 13 is prime... Well, the experiment proves the claim."

What concerns the second statement of the physicist about mathematicians, he also makes a mistake based on an inherent reductionist intentions of the majority of physicists – they try to use as little terms as possible. However, imagine if mathematicians (and physicists) instead of speaking about groups, in all cases use the term "a set with a binary operation that is associative, called multiplication, and for which an element $e$ exists, which is called the unit and which is being multiplied by any other element $a$ it yields $a$, and for any element $a$ in this set there is another element, which is called by these "pernicious" mathematicians the inverse of $a$, such that being multiplied by $a$ it yields $e$."

Another situation with the theory of non-Diophantine arithmetics was even more striking. One man said to the author:

*"I don't know your theory. I will not read anything about it. However, I completely disagree with you."*

Thus, we see that non-Diophantine arithmetics provide a mathematical base for a variety of ideas and situations in different areas, but it very difficult for people to understand even a possibility of existence for non-Diophantine arithmetics, not speaking about their essence and utilization.

It is also necessary to note that in (Burgin, 1977; 1987) non-Diophantine arithmetics are defined in a constructive way by a generative schema. Thus, it is interesting to develop logical description of these arithmetics and to give their axiomatic characterization.